\documentclass[twoside,11pt]{article}
\usepackage{graphicx,float,amscd,amsfonts,amsmath,amsthm,amssymb,latexsym,amsfonts,color}
\usepackage{epsfig,float,epstopdf,array,bm,cite,algorithmic,algorithm}
\usepackage[bookmarksnumbered,colorlinks,plainpages]{hyperref}
\usepackage{caption}
\usepackage{amsmath, latexsym}
\usepackage{cite}
\usepackage{multirow}
\newcommand{\cred}[1]{{  #1}}
\def\rank{\textrm{rank}}
 
 \newtheorem{theorem}{\bf Theorem}
 \newtheorem{lemma}{\bf Lemma}

 \newtheorem{remark}{\bf Remark}
\newtheorem{Algor}{\bf Algorithm}
\newtheorem{Example}{\bf Example}
\newcommand{\beq}{\begin{equation}}
\newcommand{\eeq}{\end{equation}}
\newcommand{\beqn}{\begin{eqnarray}}
\newcommand{\eeqn}{\end{eqnarray}}
\date{}

\begin{document}

\date{}
\title{A shift-splitting preconditioner for asymmetric saddle
point problems\thanks{The research work of the first author is supported by National
Natural Science Foundation of China (No.11961082). The work of the second author is supported by University of Guilan and the Center of Excellence for Mathematical Modelling, Optimization and Combinational Computing (MMOCC).}}
\author{Shi-Liang Wu, ~~Davod Khojasteh Salkuyeh\thanks{Corresponding
author: khojasteh@guilan.ac.ir}\\ [2mm]
{\small{\it $^{\dagger}$School of Mathematics, Yunnan Normal University, Kunming,}}\\
{\small{\it  Yunnan, 650500, P.R. China}}\\
\textit{{\small wushiliang1999@126.com}}\\[2mm]
\textit{{\small $^{\ddag}$Faculty of Mathematical Sciences, University of Guilan, Rasht, Iran }}\\
\textit{{\small $^{\ddag}$Center of Excellence for Mathematical Modelling, Optimization and Combinational}}\\
\textit{{\small Computing (MMOCC), University of Guilan, Rasht, Iran}}\\
\textit{{\small khojasteh@guilan.ac.ir}}\\[-0.5cm]}
 \maketitle
\begin{abstract}
In this paper, we execute the shift-splitting preconditioner for
asymmetric saddle point problems with its (1,2) block's
transposition unequal to its (2,1) block under the removed minus of
its (2,1) block. The proposed preconditioner is stemmed from the
shift splitting (SS) iteration method for solving asymmetric saddle
point problems, which is convergent under suitable conditions. The
relaxed version of the shift-splitting preconditioner is obtained as
well. The spectral distributions of the related preconditioned
matrices are given. Numerical experiments from the Stokes problem
are offered to show the convergence performance of these two
preconditioners.


\end{abstract}

\section{Introduction}
Nowadays, a shift splitting iteration scheme  has been successfully
used to solve the large sparse system of linear equations
\begin{equation}\label{eq:1}
 Ax=b
\end{equation}
with $A$ being non-Hermitian positive definite, which is deemed as
one of the efficient stationary solvers and is first introduced in
\cite{bai}, and works as follows:  Given an initial guess $x^{(0)}$,
for $k=0,1,2,\ldots$ until $\{x^{(k)}\}$ converges, compute
\begin{equation}\label{eq:2}
(\alpha I+ A)x^{(k+1)}=(\alpha I-A)x^{(k)}+2b,
\end{equation}
where $\alpha$ is a given positive constant. It is noteworthy that
this shift splitting iteration scheme (\ref{eq:2}) not only is
unconditionally convergent, but also can induce an economical and
effective preconditioner $P=\alpha I+ A$ for the non-Hermitian
positive definite linear system (\ref{eq:1}). This induced
precondtioner is called as the shift-splitting preconditioner. When
the shift-splitting preconditioner $P=\alpha I+ A$ together with
Krylov subspace methods are employed to solve the non-Hermitian
positive definite linear system (\ref{eq:1}), its highly efficiency
has been confirmed by numerical experiments in \cite{bai}.

Since both the shift splitting iteration scheme and the
shift-splitting preconditioner are economical and effective, they
have drawn much attention. Not only that, this approach has been
successfully extended to other practical  problems, such as  the
classical saddle point problems
\begin{equation}\label{eq:13}
\left[\begin{array}{cc}
A&B^{T}\\
-B&0
\end{array}\right]
\left[\begin{array}{c}
x\\
y
\end{array}\right]=\left[\begin{array}{c}
p\\
q
\end{array}\right],
\end{equation}
where $A\in \mathbb{R}^{n\times n}$ is symmetric positive definite (SPD),
$B\in \mathbb{R}^{m\times n}$ with $\mbox{rank}(B)=m\leq n$, see
\cite{Caoy}. Whereafter, Chen and Ma in \cite{Chen} proposed the
two-parameter shift-splitting preconditioner for  saddle point
problems (\ref{eq:13}). Based on the work in \cite{Chen}, Salkuyeh
\emph{et al.} in \cite{Salkuyeh} use the two-parameter
shift-splitting preconditioner for the saddle point problems
(\ref{eq:13}) with symmetric positive semidefinite (2, 2)-block, 
\cred{and for the same problem when the symmetry of the (1,1)-block is omitted in \cite{SalkuyehYazd}}, 
 Cao
\emph{et al.} in \cite{Caoy2} considered  the saddle point problems
(\ref{eq:13}) with nonsymmetric positive definite (1, 1)-block, Cao
and Miao in \cite{Caoy3} considered the singular nonsymmetric saddle
point problems (\ref{eq:13}), and so on.

On the other hand, combining the shift splitting technique with  the
matrix splitting technique, some new efficient preconditioners have
been developed, such as the modified shift-splitting
preconditioner\cite{Zhou}, the generalized modified shift-splitting
preconditioner\cite{Huang}, the extended shift-splitting
preconditioner \cite{Zheng}, a general class of shift-splitting
preconditioner \cite{Caoy4}, the modified generalized
shift-splitting preconditioner \cite{Huang2,DKSMRCAMWA}, the generalized double
shift-splitting preconditioner \cite{Fan}, and so on.

In this paper, we consider the asymmetric saddle point problems of the form
\begin{equation}\label{eq:14}
    \mathcal{A}{\bf x}=\left[\begin{array}{cc}
        A&B^{T}\\
        -C&0
    \end{array}\right]
    \left[\begin{array}{c}
        x\\
        y
    \end{array}\right]=\left[\begin{array}{c}
    b\\
    q
\end{array}\right]={\bf f},
\end{equation}
where $A\in \mathbb{R}^{n\times n}$ is SPD,
$B, C\in \mathbb{R}^{m\times n}$, $m\leq n$. Moreover, the matrices $B$ and $C$ are of full rank. In \cite{CaoNLAA}, Cao proposed the
augmentation block triangular preconditioner
\begin{equation}\label{PrAug}
P_{Aug}=\left[\begin{array}{cc}
A+B^TW^{-1}C&B^{T}\\
    0       & W
\end{array}\right],
\end{equation}
for the system obtaining from multiplying the second block row of  (\ref{eq:14}) by $-1$,  where $W\in\Bbb{R}^{m\times n}$ is nonsingular and such that $A+B^TW^{-1}C$ is invertible. The performance of  the preconditioner $P_{Aug}$
was compared with several preconditioners presented in \cite{CaoIJCM,CaoAPNUM,Murphy}. In \cite{LiHuangLi}, Li {\it et al.} presented the partial positive semidefinite and skew-Hermitian splitting (for short, PPSS) iteration method for the system (\ref{eq:14}). The PPSS iteration method induces the preconditioner
\begin{equation}\label{PrPPSS}
P_{PPSS}=\frac{1}{2\alpha} (\alpha I+H)(\alpha I+S),
\end{equation}
where $\alpha>0$,
\[
H=\left[\begin{array}{cc}
A & 0\\
0 & 0
\end{array}\right]\quad {\rm and} \quad
S=\left[\begin{array}{cc}
0 & B^T\\
-C & 0
\end{array}\right].
\]
Numerical results presented in \cite{LiHuangLi} show that the $P_{PPSS}$ preconditioner outperforms the classical HSS preconditioner \cite{BenziSIMAX}.  Although the shift-splitting iteration scheme and the
shift-splitting preconditioner have been successfully used to solve
the classical saddle point problems, they have not been applied to
the asymmetric saddle point problems (\ref{eq:14}). Based on this, our goal of this paper is to use the
shift-splitting iteration scheme and the shift-splitting
preconditioner for the asymmetric saddle point problems. One can see
\cite{Elman,greif1, greif2,cafieri,r,b1} for more details.
Theoretical analysis shows that the  shift-splitting iteration
method is convergent under suitable conditions and the spectral
distributions of the corresponding preconditioned matrices are
better clustered. Numerical experiments arising from a model Stokes
problem are provided to show the effectiveness of the proposed two
preconditioners.

We use the following notations throughout the paper. For a given matrix $S$, ${\cal N}({S})$ stands for the null space of $S$. The spectral radius of a square matrix $G$ is denoted by $\rho(G)$. For a vector $x\in\Bbb{C}^n$,
$x^*$ is used for the conjugate transpose of $x$. The real and imaginary parts of any $y\in\Bbb{C}$ are denoted by
$\Re(y)$ and $\Im(y)$, respectively. For two vectors $x$ and $y$, the \textsc{Matlab} notation $[x;y]$ is used for
$[x^T,y^T]^T$. Finally, for two vectors $x,y\in\Bbb{C}^n$, the standard inner product of $x$ and $y$ is denoted  by $\langle x,y \rangle=y^*x$.

The layout of this paper is organized as follows. In Section \ref{Sec2}, the
shift splitting iteration scheme and the related  shift-splitting
preconditioner are presented for the asymmetric saddle point
problems (\ref{eq:14}).  In Section \ref{Sec3}, numerical experiments are provided to
examine the convergence behaviors of the  shift-splitting
preconditioner and its relaxed version for solving the asymmetric
saddle point problems (\ref{eq:14}). Finally, some conclusions are
described in Section \ref{Sec4}.

\section{The shift-splitting method}\label{Sec2}

Here, three lemmas are given for later discussion.

\begin{lemma}\label{ExistSol} \emph{\cite{CaoNLAA}}
    The saddle point matrix
    \begin{equation}\label{CalA}
    \mathcal{A}=\left[\begin{array}{cc}
    A&B^{T}\\
    -C&0
    \end{array}\right]
    \end{equation}
    is nonsingular if and only if  $rank(B)=rank(C)=m$,
    ${\cal N}(A) \cap {\cal N}(C)=\{0\}$ and
    ${\cal N}(A^T) \cap {\cal N}(B)=\{0\}$.
\end{lemma}

\begin{lemma}\label{Root1}  \emph{\cite{Wu2}}
    Let $\lambda$ be any root of the quadratic equation $x^{2}-ax+b= 0$, where $a,b\in \mathbb{R}$. Then, $|\lambda|<1$ if and only if $|b| < 1$ and $|a| < 1+b$.
\end{lemma}

\begin{lemma}\label{Root2}\emph{\cite{Bai2}}
Let $\lambda$ be any root of the quadratic equation $x^{2 }-\phi x +
\psi = 0$, where $\phi,\psi\in \mathbb{C}$. Then, $|\lambda|<1$ if
and only if  $|\psi| < 1$ and $|\phi-\phi^{\ast}\psi| + |\psi|^{2} <
1$.
\end{lemma}

First, to guarantee the  unique solution of  the asymmetric saddle
point problems (\ref{eq:14}), Lemma \ref{NonSinA} is obtained.

\begin{lemma}\label{NonSinA}
Let $A$ be a SPD matrix and $\rank(B)=rank(C)=m$. Then, saddle point matrix $(\ref{CalA})$  is nonsingular.
\begin{proof}
It is an immediate result of Lemma \ref{ExistSol}.
\end{proof}
\end{lemma}

Similarly, for every $\alpha>0$ and under the conditions of Lemma \ref{NonSinA}, the matrix
\[
\alpha I+\mathcal{A}=\left[\begin{array}{cc}
   \alpha I+A & B^{T} \\
    -C& \alpha I \\
  \end{array}\right]
\]
is nonsingular.

Next, under the condition of Lemma \ref{NonSinA}, we can establish the shift-splitting (SS) iteration method
for solving the asymmetric saddle point problems (\ref{eq:14}). To this end, the shift-splitting of
the coefficient matrix $\mathcal{A}$ in (\ref{eq:14})  can be
constructed as follows
\begin{align*}
\mathcal{A}&=\frac{1}{2}(\alpha I+\mathcal{A})-\frac{1}{2}(\alpha
I-\mathcal{A})\\
&= \frac{1}{2}\left[\begin{array}{cc}
   \alpha I+A & B^{T} \\
    -C& \alpha I \\
  \end{array}\right]-\frac{1}{2}\left[\begin{array}{cc}
   \alpha I-A & -B^{T} \\
    C& \alpha I \\
  \end{array}\right],
\end{align*}
where $\alpha>0$ and $I$ is the identity matrix. This matrix
splitting naturally leads to the shift splitting (SS) iteration
method for solving the asymmetric saddle point problems
(\ref{eq:14}) and works as follows.

\medskip

\noindent{\it {\bf The SS iteration method}: Let the initial vector ${\bf x}^{(0)}\in
\mathbb{R}^{n+m}$ and $\alpha>0$. For
$k=0, 1, 2,\ldots$ until the iteration sequence
$\{{\bf x}^{(k)}\}_{k=0}^{+\infty}$ is converged, compute ${\bf x}^{(k+1)}$, by solving the linear system
\begin{equation}\label{eq:21}
\left[\begin{array}{cc}
   \alpha I+A & B^{T} \\
    -C& \alpha I \\
  \end{array}\right]
  {\bf x}^{(k+1)}
   =\left[\begin{array}{cc}
   \alpha I-A & -B^{T} \\
    C& \alpha I \\
  \end{array}\right]{\bf x}^{(k)}+2\left[
  \begin{array}{c}
    b \\
    q\\
  \end{array}
\right].
\end{equation}}

Clearly, the iteration matrix $M_{\alpha}$ of the SS method is
\begin{equation}\label{eq:23}
M_{\alpha}=\left[\begin{array}{cc}
   \alpha I+A & B^{T} \\
    -C& \alpha I \\
  \end{array}\right]^{-1}\left[\begin{array}{cc}
   \alpha I-A & -B^{T} \\
    C& \alpha I \\
  \end{array}\right].
\end{equation}

To study the convergence property of the SS method, the value of the
spectral radius $\rho(M_{\alpha})$ of the corresponding iteration
matrix $M_{\alpha}$ is necessary to be estimated. As is known, when
$\rho(M_{\alpha})<1$, the SS iteration method is convergent.
Thereupon, we assume that $\lambda$ is an eigenvalue of the matrix
$M_{\alpha}$ and its corresponding eigenvector is ${\bf x}=[x;y]$.
Therefore, we have
\[
\left[\begin{array}{cc}
   \alpha I-A & -B^{T} \\
    C& \alpha I \\
  \end{array}\right]\left[
  \begin{array}{c}
    x \\
    y\\
  \end{array}
\right]=\lambda\left[\begin{array}{cc}
   \alpha I+A & B^{T} \\
    -C& \alpha I \\
  \end{array}\right]\left[
  \begin{array}{c}
    x \\
    y\\
  \end{array}
\right],
\]
which is equivalent to
\begin{align}
&(\lambda-1)\alpha x+(\lambda+1)Ax+(\lambda+1)B^{T}y=0, \label{eq:23}\\
&(1+\lambda)Cx-\alpha(\lambda-1) y=0. \label{eq:24}
\end{align}

To obtain the convergence conditions of the SS method, the following
lemmas are given.

\begin{lemma}\label{Lemma5} Let the matrix $A$ be SPD and $\rank(B)=rank(C)=m$. If $\lambda$
is an eigenvalue of the matrix $M_{\alpha}$, then $\lambda\neq\pm1$.
\begin{proof}
If $\lambda=1$, then based on Eqs. (\ref{eq:23}) and
(\ref{eq:24}) we have
\begin{equation}\label{eq:25}
\left\{ \begin{aligned} &Ax+B^{T}y=0,\\
& -Cx=0.
\end{aligned} \right.
\end{equation}
Based on Lemma \ref{NonSinA}, we deduce that  $x=0$ and $y=0$.  This is a contradiction, because ${\bf x}=[x;y]=0$ can not be an
eigenvector of $M_{\alpha}$. Hence,  $\lambda\neq1$.

When $\lambda=-1$,  based on Eqs. (\ref{eq:23}) and (\ref{eq:24}) we
have  $\alpha x=0$ and  $\alpha y=0$. Since $\alpha>0$, we get $y=0$ and $x=0$, which is a contradiction, since $[x;y]$ is an eigenvector. Hence $\lambda\neq-1$.
\end{proof}
\end{lemma}

Based on the above discussion, the results in Lemma \ref{Lemma6} are right.

\begin{lemma}\label{Lemma6}  Let the conditions of Lemma \ref{Lemma5} be satisfied.
Let also $\lambda$ be an eigenvalue of $M_{\alpha}$ and ${\bf x}=[x;y]$ be
the corresponding eigenvector. Then $x \neq0$. Moreover, if $y = 0$,
then $|\lambda| < 1$.
\begin{proof}
When $x=0$, from (\ref{eq:24}) we have
$\alpha(\lambda-1) y=0$. Based on Lemma \ref{Lemma5}, $\lambda\neq 1$. Therefore, $y=0$. This
contradicts with the nonzero eigenvector ${\bf x}=[x;y]$. Hence $x
\neq0$.

When $y=0$, based on Eq. (\ref{eq:23}) we get
\[
( \alpha I+A)^{-1}(\alpha I-A)x=\lambda x.
\]
Therefore, using the Kellogg's lemma (see \cite[page 13]{Marchuk}) we deduce
\[
|\lambda|\leq \|( \alpha I+A)^{-1}(\alpha I-A)\|_{2}<1,
\]
which completes the proof.
\end{proof}
\end{lemma}

For later use we define the set ${\cal S}$ as
\[
{\cal S}=\{x\in\Bbb{C}^n: {\bf x}=[x;y] \text{ is an eigenvector of } M_{\alpha} \text{ with } \|x\|_2=1\}.
\]
It follows from Lemma \ref{Lemma6} that the members of ${\cal S}$ are nonzero.

\begin{theorem}\label{Thm1}
Let the conditions of Lemma \ref{Lemma5} be satisfied.
For every $x\in{\cal S}$, let $a(x)=x^{\ast}Ax$, $s(x)=\Re(x^HB^TCx)$ and $t(x)=\Im(x^HB^TCx)$.
For each $x\in{\cal S}$, if $s(x)>0$ and
\begin{equation}\label{EqCond}
|t(x)|<a(x)\sqrt{s(x)},
\end{equation}
then
\begin{equation*}
\rho(M_{\alpha})<1,\quad \forall \alpha>0,
\end{equation*}
which implies that the SS iteration method $(\ref{eq:21})$ converges
to the unique solution of the asymmetric saddle point problems
$(\ref{eq:14})$.
\begin{proof}
Based on Lemma \ref{Lemma5}, from (\ref{eq:24}) we have
\begin{equation}\label{eq:27}
y=\frac{\lambda+1}{\alpha(\lambda-1)}Cx.
\end{equation}
Substituting (\ref{eq:27}) into (\ref{eq:23}) leads to
\begin{equation}\label{eq:28}
(\lambda-1)\alpha
x+(\lambda+1)Ax+\frac{(\lambda+1)^{2}}{\alpha(\lambda-1)}B^{T}Cx=0.
\end{equation}
Let $\|x\|_{2}=1$. Pre-multiplying $x^{\ast}$ to the both sides  of
Eq. (\ref{eq:28}) leads to
\begin{equation}\label{EqQuadEq}
\alpha^{2}(\lambda-1)^{2}+\alpha(\lambda^{2}-1)x^{\ast}Ax+(\lambda+1)^{2}x^{\ast}B^{T}Cx=0,
\end{equation}
which is equivalent to
\begin{equation}\label{eq:29}
\alpha^{2}(\lambda-1)^{2}+\alpha(\lambda^{2}-1)a+(\lambda+1)^{2}(s(x)+t(x)i)=0.
\end{equation}
For the sake simplicity in notations, we use $s$, $t$ and $a$ for $s(x)$, $t(x)$ and $a(x)$, respectively.
It follows from Eq. (\ref{eq:29}), that
\begin{equation}\label{eq:210}
\lambda^{2}+\frac{2(s+ti-\alpha^{2})}{\alpha^{2}+\alpha
a+s+ti}\lambda+\frac{\alpha^{2}-\alpha a+s+ti}{\alpha^{2}+\alpha
a+s+ti}=0.
\end{equation}
Next, we will discuss two aspects: $t = 0$ and $t\neq0$.

When $t = 0$, from (\ref{eq:210}), we get
\begin{equation}\label{eq:211}
\lambda^{2}+\frac{2(s-\alpha^{2})}{\alpha^{2}+\alpha
a+s}\lambda+\frac{\alpha^{2}-\alpha a+s}{\alpha^{2}+\alpha a+s}=0.
\end{equation}
By simples computations, we have
\begin{equation}\label{eq:212}
\Big|\frac{\alpha^{2}-\alpha a+s}{\alpha^{2}+\alpha a+s}\Big|<1
\end{equation}
and
\begin{equation}\label{eq:213}
\Big|\frac{2(s-\alpha^{2})}{\alpha^{2}+\alpha a+s}\Big|<1+
\frac{\alpha^{2}-\alpha a+s}{\alpha^{2}+\alpha a+s}.
\end{equation}
Based on Lemma \ref{Root1}, the inequalities (\ref{eq:212}) and (\ref{eq:213})
imply that the roots of the real quadratic equation (\ref{eq:211})
satisfy $|\lambda| < 1$.

If $t\neq0$, then Eq. (\ref{eq:210}) can be written as $\lambda^2+\phi \lambda +\psi=0$, where
\[
\phi=\frac{2(s+ti-\alpha^{2})}{\alpha^{2}+\alpha a+s+ti}\
\mbox{and}\ \psi=\frac{\alpha^{2}-\alpha a+s+ti}{\alpha^{2}+\alpha
a+s+ti}.
\]
By some calculations, we get
\begin{align*}
\phi-\phi^{\ast}\psi&=\frac{2(s+ti-\alpha^{2})}{\alpha^{2}+\alpha
a+s+ti}-\frac{2(s-ti-\alpha^{2})}{\alpha^{2}+\alpha
a+s-ti}\cdot\frac{\alpha^{2}-\alpha a+s+ti}{\alpha^{2}+\alpha
a+s+ti}\\
&=\frac{2(s-\alpha^{2}+ti)}{\alpha^{2}+\alpha
a+s+ti}\cdot\frac{\alpha^{2}+\alpha a+s-ti}{\alpha^{2}+\alpha
a+s-ti}-\frac{2(s-\alpha^{2}-ti)}{\alpha^{2}+\alpha
a+s-ti}\cdot\frac{\alpha^{2}-\alpha a+s+ti}{\alpha^{2}+\alpha a+s+ti}\\
&=2\Big[\frac{(s-\alpha^{2}+ti)(\alpha^{2}+\alpha
a+s-ti)}{(\alpha^{2}+\alpha
a+s)^{2}+t^{2}}+\frac{(\alpha^{2}-s+ti)(\alpha^{2}-\alpha
a+s+ti)}{(\alpha^{2}+\alpha a+s)^{2}+t^{2}}\Big]\\
&=2\frac{(s-\alpha^{2}+ti)(\alpha^{2}+\alpha
a+s-ti)+(\alpha^{2}-s+ti)(\alpha^{2}-\alpha
a+s+ti)}{(\alpha^{2}+\alpha
a+s)^{2}+t^{2}}\\
&=4\frac{\alpha a(s-\alpha^{2})+2\alpha^{2} ti}{(\alpha^{2}+\alpha
a+s)^{2}+t^{2}}.
\end{align*}
Further, we have
\begin{equation} \label{eq:217}
\begin{aligned}
&|\psi|=\sqrt{\frac{(\alpha^{2}-\alpha
a+s)^{2}+t^{2}}{(\alpha^{2}+\alpha a+s)^{2}+t^{2}}}<1,\\
&|\phi-\phi^{\ast}\psi|=\frac{4\sqrt{\alpha^2
a^2(s-\alpha^{2})^{2}+4t^{2}\alpha^{4}}}{(\alpha^{2}+\alpha
a+s)^{2}+t^{2}}.
\end{aligned}
\end{equation}
Based on Lemma \ref{Root2}, the necessary and sufficient condition for
$|\lambda|<1$ is
\begin{equation}\label{eq:218}
|\phi-\phi^{\ast}\psi|+|\psi|^{2}<1.
\end{equation}
Substituting (\ref{eq:217}) into (\ref{eq:218}) and solving the
inequality (\ref{eq:218})  for $t$, gives
$
|t|< a\sqrt{s},
$
which completes the proof.
\end{proof}
\end{theorem}

According to the definition of $t(x)$, we have
\begin{eqnarray*}
|t(x)|&=&|\Im(x^HB^TCx)|=|\langle B^TCx,x \rangle| \\
      & \leq& \|B^TCx\|_2 \|x\|_2  \qquad\qquad\qquad\qquad\qquad \text{(Cauchy-Schwarz inequality)} \\
      & \leq& \|B^TC\|_2 \|x\|_2=\|B^TC\|_2.
\end{eqnarray*}
Also we have $a(x)=x^{\ast}Ax\geq \lambda_{\min}(A)$, where $\lambda_{\min}(A)$ is the smallest eigenvalue of $A$. Therefore, the inequality (\ref{EqCond}) can be replaced by
\[
\|B^TC\|_2 \leq \lambda_{\min}(A) \sqrt{s(x)}.
\]

In the special case that $C=kB$ with $k>0$, we can state the following theorem.

\begin{theorem}\label{Thm2}
Let the conditions of Lemma \ref{Lemma5} be satisfied  and $C=kB$ with $k>0$.
Then $\rho(M_{\alpha})<1$,  $\forall \alpha>0$, which implies that the SS iteration method $(\ref{eq:21})$ converges
to the unique solution of the asymmetric saddle point problems
$(\ref{eq:14})$.
\begin{proof}
If $C=kB$ with $k>0$, then the matrix $B^TC=kB^TB$ is symmetric positive semidefinite. Therefore, we have $s(x)=kx^*B^TBx\geq 0$ and $t(x)=0$. According to Theorem \ref{Thm1}, all we need is to prove the convergence for the case that $s(x)=0$. If $s(x)=0$, then  we get $Bx=0$. Now, from Eq. (\ref{CalA}) we deduce that
\[
\alpha^{2}(\lambda-1)^{2}+\alpha(\lambda^{2}-1)a=0,
\]
which is equivalent to
\[
\alpha(\lambda-1) ( \alpha^{2}(\lambda-1)+\alpha(\lambda+1)a)=0.
\]
Now, since $\alpha> 0$ and $\lambda \neq 1$ (from Lemma \ref{Lemma5}), we deduce that
\[
\alpha^{2}(\lambda-1)+\alpha(\lambda+1)a=0,
\]
which gives the following equation for $\lambda$
\[
\lambda=\frac{\alpha-a}{\alpha+a}.
\]
Therefore, since $a=x^*Ax>0$, we conclude that $|\lambda|<1$, which completes the proof.
\end{proof}
\end{theorem}

\begin{remark}
When $k=1$, Theorem \ref{Thm2} is the main result in \cite{Caoy}. That is to say, Theorems \ref{Thm1} and \ref{Thm2}
are generalizations Theorem 2.1 in \cite{Caoy}.
\end{remark}

Finally, we consider the preconditioner induced by the SS iteration method $(\ref{eq:21})$.
As is known, the advantage of matrix splitting technique often is twofold: one is to result in a
splitting iteration method and the other is to induce a splitting
preconditioner for improving the convergence speed of Krylov
subspace methods in \cite{bai}. Based on the SS iteration method
$(\ref{eq:21})$, the corresponding shift-splitting preconditioner
can be defined by
\[
P_{SS}=\frac{1}{2}\left[\begin{array}{cc}
   \alpha I+A & B^{T} \\
    -C& \alpha I \\
  \end{array}\right].
\]
Since the multiplicative factor $\frac{1}{2}$ in the preconditioner
$P_{SS}$ has no effect and can be removed when $P_{SS}$ is used as a
preconditioner, in the implementations, we only consider the
shift-splitting preconditioner $P_{SS}$ without  the multiplicative
factor $\frac{1}{2}$. In this case, using $P_{SS}$ with Krylov
subspace methods (such as GMRES, or its restarted version GMRES($k$)),
a vector of the form
\[
z=P_{SS}^{-1}r
\]
needs to be computed.

Let $z=[z_{1};z_{2}]$ and $r=[r_{1};r_{2}]$. Then $z=P_{SS}^{-1}r$ is equal to
\begin{equation}\label{eq:31}
\left[
  \begin{array}{c}
    z_{1} \\
    z_{2}\\
  \end{array}
\right]=\left[\begin{array}{cc}
   I & 0 \\
   \frac{1}{\alpha}C& I \\
  \end{array}\right]\left[\begin{array}{cc}
   \alpha I+A+ \frac{1}{\alpha}B^{T}C& 0\\
    0& \alpha I\\
  \end{array}\right]^{-1}\left[\begin{array}{cc}
    I & -\frac{1}{\alpha}B^{T} \\
   0&  I \\
  \end{array}\right]\left[
  \begin{array}{c}
    r_{1} \\
    r_{2}\\
  \end{array}
\right].
\end{equation}
Based on Eq. (\ref{eq:31}), the following algorithm can be used  to
obtain the vector $z$.

\medskip

\begin{Algor}\label{Algor1}\rm
Let $z=[z_{1};z_{2}]$ and
$r=[r_{1};r_{2}]$. Compute $z$ by the following procedure\\
1. Compute $t=r_{1}-\frac{1}{\alpha}B^{T}r_{2}$;\\
2. Solve  $(\alpha I+A+ \frac{1}{\alpha}B^{T}C)z_{1}=t$ for $z_1$;\\
3. Compute $z_{2}=\frac{1}{\alpha}(Cz_{1}+r_{2})$.
\end{Algor}

In Step 2 of Alg. \ref{Algor1}, in general the matrix $\alpha I+A+
\frac{1}{\alpha}B^{T}C$ is indefinite, hence the corresponding
system can be solved exactly using the LU factorization or inexactly
using a Krylov subspace method like GMRES or its restarted version.
However, when $C=kB$ with $k>0$, this matrix is of the form $\alpha
I+A+\frac{k}{\alpha}B^{T}B$ which is SPD. Therefore, the
corresponding system can be solved exactly using the Cholesky
factorization or inexactly using the conjugate gradient (CG) method.

\cred{In general the matrix $\alpha I+A+\frac{1}{\alpha}B^{T}C$  is dense (because of the term $B^TC$) and solving the corresponding system by a direct method may be impractical. Hence, it is recommended to solve the system by an iteration method, as we will shortly do in the section of the numerical experiments. From theoretical point of view, when $\alpha=0$ the preconditioner $P_{SS}=\alpha I+\mathcal{A}$ coincides with the coefficient matrix of original system. In this case, implementation of the preconditioner would be as difficult as solving the original system. Hence, it is better to choose a small value of $\alpha$ to obtain a more well-conditioned matrix. Since the condition of the matrix 
$\alpha I+A+\frac{1}{\alpha}B^{T}C$ strongly depends on the term $\frac{1}{\alpha}B^{T}C$, similar to \cite{Caoy2,Golub} we choose 
the parameter $\alpha$ equals to 
\[
\alpha_{est}=\frac{\|B^TC\|_2}{\|A\|_2},
\]
which balances the matrices $A$ and $B^TC$.}

When Krylov subspace methods together with the preconditioner
$P_{SS}=\alpha I+\mathcal{A}$ are applied  to solve the asymmetric
saddle point problems $(\ref{eq:14})$, \cred{we need to establish} the
spectral distribution of the preconditioned matrix $P_{SS}^{-1}A$ to
investigate the convergence performance of the preconditioner
$P_{SS}$ for Krylov subspace methods.

The following theorem on the spectral distribution of the
preconditioned matrix $P_{SS}^{-1}\mathcal{A}$ can be obtained.
\begin{theorem}\label{Thm3} Let the conditions of
Theorem \ref{Thm1} or \ref{Thm2}  be satisfied. Then the
preconditioned matrix $P_{SS}^{-1}\mathcal{A}$ are positive stable
for $\alpha>0$ and \cred{its the eigenvalues} satisfy $|\lambda|<1$,
where $\lambda$ denotes  the eigenvalue of the preconditioned matrix
$P_{SS}^{-1}\mathcal{A}$.
\begin{proof}
\cred{It follows from
\begin{align*}
2P_{SS}^{-1}\mathcal{A}=I-M_{\alpha},
\end{align*}
that for each $\mu \in \sigma(M_{\alpha})$, there is a $\lambda\in \sigma(P_{SS}^{-1}\mathcal{A})$, such that $2\lambda=1-\mu$. Therefore, we
\begin{eqnarray*}
\frac{\mu}{2} &=& \frac{1}{2} - \lambda = \frac{1}{2} - \Re(\lambda) - i \Im(\lambda).  
\end{eqnarray*}
Hence, from the fact that $|\mu|<1$ we conclude
\[
(\frac{1}{2} - \Re(\lambda))^2 +(\Im(\lambda))^2< \frac{1}{4},
\]
which shows that the eigenvalues of the preconditioned matrix $P_{SS}^{-1}\mathcal{A}$ are contained in a circle with radius $\frac{1}{2}$ centered at $(\frac{1}{2},0)$. Hence,  the real parts of the eigenvalues of the matrix
$P_{SS}^{-1}\mathcal{A}$ are all positive. This means that the matrix $P_{SS}^{-1}\mathcal{A}$ is
positive stable for $\alpha>0$. On the other hand, from $|\mu|< 1$  we deduce that
\[
2|\lambda|=|1-\mu|\leq 1+|\mu|< 2,
\]
which completes the proof.}
\end{proof}
\end{theorem}

Here, we present a relaxed version of the shift-splitting
preconditioner as well, which is defined by
\[
P_{RSS}=\left[\begin{array}{cc}
   A & B^{T} \\
    -C& \alpha I \\
  \end{array}\right].
\]
Similarly, using $P_{RSS}$ with Krylov  subspace methods (such as
GMRES, or its restarted version GMRES($k$)), a vector of the form
\[
z=P_{RSS}^{-1}r
\]
has to be computed as well.  Let $z=[z_{1};z_{2}]$ and
$r=[r_{1};r_{2}]$. Then,  we have
\begin{equation}
\left[
  \begin{array}{c}
    z_{1} \\
    z_{2}\\
  \end{array}
\right]=\left[\begin{array}{cc}
   I & 0 \\
   \frac{1}{\alpha}C& I \\
  \end{array}\right]\left[\begin{array}{cc}
   A+ \frac{1}{\alpha}B^{T}C& 0\\
    0& \alpha I\\
  \end{array}\right]^{-1}\left[\begin{array}{cc}
    I & -\frac{1}{\alpha}B^{T} \\
   0&  I \\
  \end{array}\right]\left[
  \begin{array}{c}
    r_{1} \\
    r_{2}\\
  \end{array}
\right].
\end{equation}

Based on Alg. \ref{Algor1}, by a simple modification,  we obtain Alg. \ref{Algor2} to obtain the vector $z$ as follows.
\begin{Algor}\label{Algor2}\rm
Let $z=[z_{1};z_{2}]$ and $r=[r_{1};r_{2}]$. Compute $z$ by the following procedure\\
1. Compute $t=r_{1}-\frac{1}{\alpha}B^{T}r_{2}$;\\
2. Solve  $(A+ \frac{1}{\alpha}B^{T}C)z_{1}=t$;\\
3. Compute $z_{2}=\frac{1}{\alpha}(Cz_{1}+r_{2})$.
\end{Algor}

In the same way, we can obtain the spectral distribution of the
preconditioned matrix $P_{RSS}^{-1}\mathcal{A}$, as follows.
\begin{theorem}\label{Thm4}
Let the conditions of Theorem \ref{Thm1}  be satisfied. Then the
preconditioned matrix $P_{RSS}^{-1}\mathcal{A}$ has an eigenvalue
$1$ with algebraic multiplicity $n$ and the remaining eigenvalues
are the eigenvalues of matrix $\frac{1}{\alpha}C(A+
\frac{1}{\alpha}B^{T}C)^{-1}B^{T}$.
\begin{proof}
By calculation, we get
\begin{align*}
P_{RSS}^{-1}\mathcal{A}=&\left[\begin{array}{cc}
   I & 0 \\
   \frac{1}{\alpha}C& I \\
  \end{array}\right]\left[\begin{array}{cc}
   A+ \frac{1}{\alpha}B^{T}C& 0\\
    0& \alpha I\\
  \end{array}\right]^{-1}\left[\begin{array}{cc}
    I & -\frac{1}{\alpha}B^{T} \\
   0&  I \\
  \end{array}\right]\left[\begin{array}{cc}
   A & B^{T} \\
    -C& 0 \\
  \end{array}\right]\\
  =&\left[\begin{array}{cc}
   (A+ \frac{1}{\alpha}B^{T}C)^{-1} & 0 \\
   \frac{1}{\alpha}C(A+ \frac{1}{\alpha}B^{T}C)^{-1}& \frac{1}{\alpha} I\\
  \end{array}\right]\left[\begin{array}{cc}
    I & -\frac{1}{\alpha}B^{T} \\
   0&  I \\
  \end{array}\right]\left[\begin{array}{cc}
   A & B^{T} \\
    -C& 0 \\
  \end{array}\right]\\
  =&\left[\begin{array}{cc}
   (A+ \frac{1}{\alpha}B^{T}C)^{-1} &  -\frac{1}{\alpha}(A+ \frac{1}{\alpha}B^{T}C)^{-1}B^{T} \\
   \frac{1}{\alpha}C(A+ \frac{1}{\alpha}B^{T}C)^{-1}&-\frac{1}{\alpha^{2}}C(A+ \frac{1}{\alpha}B^{T}C)^{-1}B^{T}+  \frac{1}{\alpha} I\\
  \end{array}\right]\left[\begin{array}{cc}
   A & B^{T} \\
    -C& 0 \\
  \end{array}\right]\\
=&\left[\begin{array}{cc}
   I &  (A+ \frac{1}{\alpha}B^{T}C)^{-1}B^{T} \\
   0&\frac{1}{\alpha}C(A+ \frac{1}{\alpha}B^{T}C)^{-1}B^{T}\\
  \end{array}\right].
\end{align*}
Therefore,  the proof is completed.
\end{proof}
\end{theorem}

\cred{Obviously, for each $\alpha>0$ the preconditioner $P_{RSS}$ is more closer than the preconditioner $P_{SS}$ to the original matrix ${\cal A}$. However, the subsystem appeared in the implementation of  the $P_{SS}$ preconditioner in a Krylov subspace method is more well-conditioned than that of $P_{RSS}$. Hence, it is recommend to apply  the preconditioner $P_{SS}$ when the subsystems are solved inexactly using an iteration method  and the $P_{RSS}$ when the subsystems are solved exactly using direct method.}

\section{Numerical experiments}\label{Sec3}

In this section, we present some numerical experiments to demonstrate the performance of the
shift-splitting preconditioner. In the
meantime, the numerical comparison are provided to show the
advantage of the shift-splitting preconditioner ($P_{SS}$) and its relaxed version ($P_{RSS}$) \cred{over the PPSS preconditioner given by Eq. \eqref{PrPPSS} (denoted by  $P_{PPSS}$) and the augmentation block triangular preconditioner given by Eq. \eqref{PrAug} (denoted by $P_{Aug}$)}.
In our computations, we apply the flexible GMRES (FGMRES) \cite{FGMRES,SaadBook}  together with these four
preconditioners to solve the assymmetric saddle point systems
(\ref{eq:14}) and adjust the right-hand side ${\bf f}$ such that the exact
solution is a vector of all ones. The iterations start with a zero vector as an initial guess and are stopped when the numbers of
iteration exceeds 1000 or
\[
R_k=\frac{\|{\bf f}-\mathcal{A} {\bf x}^{(k)}\|_{2}}{\| {\bf f}\|_{2}}\leq 10^{-7},
\]
where ${\bf x}^{(k)}$ is the computed solution at iteration $k$. In the implementation of the preconditioners the subsystems are solved inexactly using the iterations method. When the coefficient matrix is SPD, the corresponding system is solved using the conjugate gradient (CG) method, otherwise by the restarted version of GMRES(10). For the subsystems, the iteration is stopped as soon as the residual 2-norm is reduced by a factor of $10^2$ and the maximum number of iterations is set to be 100. Similar to the outer iterations, a null vector is used as an initial guess. Finally, for the augmentation block triangular preconditioner the matrix $W$ is set to be $W=\alpha I$ with $\alpha>0$. In this case, the preconditioner $P_{Aug}$ takes the following form
\[
P_{Aug}=\left[\begin{array}{cc}
A+\frac{1}{\alpha}B^TC& B^{T}\\
0       & \alpha I
\end{array}\right],
\]
For all the methods the optimal value of parameter are obtained experimentally (denoted by $\alpha_*$) and are the ones resulting in the least numbers of iterations. \cred{We also report the numerical results for the parameter
$\alpha_{est}={\|B^TC\|_2}/{\|A\|_2}$.
}  

We present the numerical results in the tables. In the tables, ``CPU" and ``Iters" stand for the elapsed CPU time (in second) and the number of iterations for the convergence. A dagger ($\dag$) means that the iteration has not converged in 1000 iterations.
All runs are implemented in \textsc{Matlab} R2017, equipped with a Laptop with 1.80 GHz central processing unit (Intel(R) Core(TM) i7-4500), 6 GB memory and Windows 7 operating system.

\begin{Example}\label{Ex1} \rm Let the asymmetric saddle point
 problems (\ref{eq:14}) be given by
\begin{equation*}
A=\left[\begin{array}{cc}
   I\otimes T+T\otimes I & 0 \\
   0& I\otimes T+T\otimes I\\
  \end{array}\right]\in\mathbb{R}^{2s^{2}\times2s^{2}}
\end{equation*}
and
\begin{equation*}
B^{T}=\left[\begin{array}{cc}
   I\otimes F \\
   F\otimes I\\
  \end{array}\right]\in\mathbb{R}^{2s^{2}\times s^{2}},\ C=kB,
\end{equation*}
with
\[
T=\frac{\mu}{h^{2}}\mbox{tridiag}(-1,2,-1)\in\mathbb{R}^{s\times
s},\quad F=\frac{1}{h}\mbox{tridiag}(-1,1,0)\in\mathbb{R}^{s\times
s},k>0.
\]
where $\otimes$ denotes the Kronecker product and $h={1}/{(s+1)}$
is the discretization mesh-size. Therefore, the total number of variables $n=3s^{2}$.

This asymmetric saddle point problems (\ref{eq:14}) can be obtained by using the upwind scheme to
discretize the Stokes problem in the region
$\Omega=(0,1)\times(0,1)\subset \mathbb{R}^{2}$ with its boundary
being $\partial\Omega$: find $u$ and $p$ such that
\begin{equation*}
\left\{ \begin{aligned}
         -\mu\Delta u+\nabla p &=f, \ \mbox{in}\ \Omega,\\
                  \nabla\cdot u&=g,\ \mbox{in}\ \Omega,\\
                  u&=0,\ \mbox{on} \ \partial\Omega,\\
                  \int_{\Omega} p(x)dx&=0,
                          \end{aligned} \right.
\end{equation*}
where $\mu$, $\Delta$, $u$ and $p$ are the viscosity scalar, the
componentwise Laplace operator, a vector-valued function
representing the velocity, and a scalar function representing the
pressure, respectively.

We set $s=16,32,64,128,256$ and $k=2$. Generic properties of the test matrices are presented in Table \ref{Tbl1}. In this table, $nnz(.)$ stands for number of nonzero entries of the matrix. Numerical results for $\mu=1$ and $\mu=0.1$ are presented in the Tables \ref{Tbl2} and \ref{Tbl3}, respectively. From the numerical results in Tables \ref{Tbl2}-\ref{Tbl3}, it is easy to find that the computational efficiency of GMRES
can not be satisfy when it is directly used to solve the asymmetric saddle point problems
(\ref{eq:14}). Whereas, FGMRES together with these four preconditioners for solving the asymmetric saddle point problems
(\ref{eq:14}) can rapidly converge. This also confirms that all four preconditioners indeed can improve the convergence speed of
GMRES. Among the preconditioners, $P_{SS}$ and $P_{RSS}$ outperform the others from the iteration steps and the CPU time point of review. \cred{On the other hand we observe the parameter $\alpha_{est}$ often gives suitable results, especially for large problems.}   

\begin{table}
    \centering\caption{Matrix properties  for Example \ref{Ex1}.\label{Tbl1}}\vspace{-0.5cm}
    \begin{center}
        \scalebox{1.}
        {
            \begin{tabular}{|c|c|c|c|c|c|} \hline
                $s$ &  $n$  &  $m$   & $nnz(A)$ & $nnz(B)$ & $nnz(C)$ \\ \hline
                16  &   512    &  256   & 2432    & 992    & 992    \\ [0.2cm]
                32  &       2048   &  1024  & 9984    & 4032   & 4032   \\ [0.2cm]
                64  &       8192   &  4096  & 40448   & 16256  & 16256  \\ [0.2cm]
                128  &      32768  &  16384 & 162816  & 65280  & 65280  \\ [0.2cm]
                256  &      131072 &  65536 & 653312  & 261632 & 261632 \\ \hline
            \end{tabular}
        }
    \end{center}

\end{table}
%----------------------------------------------------------------------------------------------
\begin{table} \centering
    \centering\caption{Numerical results of FGMRES for Example \ref{Ex1} with $\mu=1$.\label{Tbl2}}
    \scalebox{0.8}{
        \begin{tabular} {|l||l|l|l|l|l|l||l|l|l|l|} \hline
            $s$   &      ~~~~~~~~~      &No Prec. &$P_{SS}$~~~~~~ & $P_{RSS}$ ~~ & $P_{PPSS}$ ~~ &$P_{Aug}$~~ &&\cred{$P_{SS}$}~~~~~~& \cred{$P_{RSS}$} ~~\\ \hline
            16    & $\alpha_*$  &    --     &  0.10   & 0.20    &  98.50  &   0.11  &  $\alpha_{est}$  & 2.03   & 2.03\\
                  &   Iters     &   133     &  8      & 8       &  38     &   21    &    Iters         & 12     & 11 \\
                  &   CPU       &   0.13    &  0.03   & 0.02    &  0.05   &   0.06  &    CPU           & 0.03   & 0.03 \\
                  &   $R_k$     &   8.1e-8  & 8.4e-8  & 5.9e-8  & 1.0e-7  &  7.5e-8 &    $R_k$         & 6.6e-8 & 9.3e-8\\ \hline
                                                                                                  
            32    & $\alpha_*$  &    --     &  0.20   & 0.34    &  100.60 &   0.10  &  $\alpha_{est}$  & 2.01   & 2.01 \\
                  &   Iters     &   285     &  9      & 9       &  45     &   21    &    Iters         & 13     & 12 \\
                  &   CPU       &   2.93    &  0.06   & 0.05    &  0.14   &   0.14  &    CPU           & 0.06   & 0.06 \\
                  &   $R_k$     &   9.6e-8  & 2.4e-8  & 9.6e-8  & 1.0e-7  &  8.9e-8 &    $R_k$         & 5.2e-8 & 7.4e-8\\ \hline
                                                                                                  
            64    & $\alpha_*$  &    --     &  0.60   & 1.50    &  102.20 &   0.37  &  $\alpha_{est}$  & 2.01   & 2.01 \\
                  &   Iters     &   617     &  12     & 12      &  63     &   29    &    Iters         & 14     & 13 \\
                  &   CPU       &   36.20   &  0.39   & 0.3     &  0.89   &   0.74  &    CPU           & 0.38   & 0.36\\
                  &   $R_k$     &   9.6e-8  & 7.5e-8  & 8.2e-8  & 9.3e-8  &  8.2e-8 &    $R_k$         & 5.6e-8 & 6.4e-8\\ \hline
                                                                                                  
            128   & $\alpha_*$  &    --     &  0.60   & 0.64    &  103.90 &   4.20  &  $\alpha_{est}$  & 2.02   & 2.02\\
                  &   Iters     &   $\dag$  &  22     & 23      &  111    &   31    &    Iters         & 24     & 23 \\
                  &   CPU       &   --      &  2.37   & 2.48    &  8.69   &   3.19  &    CPU           & 2.55   & 2.33\\
                  &   $R_k$     &   --      & 8.4e-8  & 8.5e-8  & 8.8e-8  &  6.5e-8 &    $R_k$         & 5.2e-8 & 5.4e-8\\ \hline
                                                                                                  
            256   & $\alpha_*$  &    --     &  1.39   & 1.39    &  102.00 &   22.00 &  $\alpha_{est}$  & 2.02   & 2.02\\
                  &   Iters     &  $\dag$   &  57     & 52      &  217    &   78    &    Iters         & 64     & 54 \\
                  &   CPU       &   --      &  34.89  & 32.38   &  175.49 &   47.18 &    CPU           & 40.66  & 33.84 \\
                  &   $R_k$     &   --      & 9.5e-8  & 8.1e-8  & 9.9e-8  &  7.1e-8 &    $R_k$         & 9.5e-8 & 4.2e-8 \\ \hline

        \end{tabular}}
    \end{table}

    %-----------------------------------------------------------------------------------------------------------
    \begin{table} \centering
        \centering\caption{Numerical results of FGMRES for Example \ref{Ex1} with $\mu=0.1$.\label{Tbl3}}
            \scalebox{0.8}{
      \begin{tabular} {|l||l|l|l|l|l|l||l|l|l|l|} \hline
            $s$   &      ~~~~~~~~~      &No Prec. &$P_{SS}$~~~~~~ & $P_{RSS}$ ~~ & $P_{PPSS}$ ~~ &$P_{Aug}$~~ && \cred{$P_{SS}$}~~~~~~& \cred{$P_{RSS}$} ~~\\ \hline
                  16    & $\alpha_*$  &    --   &  0.25   & 0.25   &  15.40  &   0.53  &  $\alpha_{est}$  & 18.34  & 18.34\\          
                      &   Iters     &   117     &  8      & 8      &  36     &   17    &    Iters         & 28     & 12 \\           
                      &   CPU       &   0.16    &  0.02   & 0.02   &  0.04   &   0.04  &    CPU           & 0.02   & 0.02 \\         
                      &   $R_k$     &   8.9e-8  & 1.5e-8  & 1.4e-8 & 9.4e-8  &  8.6e-8 &    $R_k$         & 8.2e-8 & 6.6e-8\\ \hline 
                                                                                                                                     
                32    & $\alpha_*$  &    --     &  0.23   & 0.23   &  29.80  &   2.42  &  $\alpha_{est}$  & 19.45  & 19.45 \\         
                      &   Iters     &   238     &  11     & 11     &  56     &   20    &    Iters         & 31     & 13 \\           
                      &   CPU       &   1.67    &  0.07   & 0.07   &  0.16   &   0.11  &    CPU           & 0.08   & 0.06 \\         
                      &   $R_k$     &   9.0e-8  & 9.0e-8  & 5.6e-8 & 9.2e-8  &  1.0e-7 &    $R_k$         & 6.9e-8 & 5.3e-8\\ \hline 
                                                                                                                                     
                64    & $\alpha_*$  &    --     &  1.50   & 2.1    &  53.20  &   4.60  &  $\alpha_{est}$  & 19.87  & 19.87 \\         
                      &   Iters     &   483     &  11     & 11     &  86     &   26    &    Iters         & 32     & 14 \\           
                      &   CPU       &   22.48   &  0.26   & 0.25   &  0.95   &   0.65  &    CPU           & 0.46   & 0.38\\          
                      &   $R_k$     &   9.9e-8  & 9.7e-8  & 5.8e-8 & 9.6e-8  &  9.4e-8 &    $R_k$         & 8.8e-8 & 4.2e-8\\ \hline 
                                                                                                                                     
                128   & $\alpha_*$  &    --     &  4.90   & 6.4    &  92.80  &   19.10 &  $\alpha_{est}$  & 19.98  & 19.98\\          
                      &   Iters     &   908     &  18     & 19     &  129    &   39    &    Iters         & 33     & 20 \\           
                      &   CPU       &   302.64  &  1.91   & 1.96   &  7.07   &   4.06  &    CPU           & 3.07   & 2.17\\          
                      &   $R_k$     &   9.9e-8  &  9.2e-8 & 7.2e-8 & 9.9e-8  &  9.9e-8 &    $R_k$         & 7.4e-8 & 9.2e-8\\ \hline 
                                                                                                                                     
                256   & $\alpha_*$  &    --     &  10.90  & 12.96  &  131.00 &   25.90 &  $\alpha_{est}$  & 20.05  & 20.05\\          
                      &   Iters     &  $\dag$   &  30     & 37     &  192    &   90    &    Iters         & 37     & 46 \\           
                      &   CPU       &   --      &  26.03  & 22.73  &  151.26 &   55.38 &    CPU           & 23.26  & 29.10 \\        
                      &   $R_k$     &   --      & 9.0e-8  & 9.6e-8 & 9.7e-8  &  7.8e-8 &    $R_k$         & 6.2e-8 & 9.1e-8 \\ \hline

            \end{tabular}}
        \end{table}

In the sequel, we investigate the spectral
distribution of four preconditioned matrices
$P_{SS}^{-1}\mathcal{A}$, $P_{RSS}^{-1}\mathcal{A}$, $P_{PPSS}^{-1}\mathcal{A}$ and
$P_{Aug}^{-1}\mathcal{A}$. To do so, we set $s=16$  and use the optimal value of the parameters given in Tables \ref{Tbl2} and \ref{Tbl3}.
Figs. \ref{Fig1}-\ref{Fig2} plot the spectral distribution of the
matrices. Fig. \ref{Fig1} plots the spectral distribution of five
matrices $\mathcal{A}$, $P_{SS}^{-1}\mathcal{A}$, $P_{RSS}^{-1}\mathcal{A}$
$P_{PPSS}^{-1}\mathcal{A}$ and $P_{Aug}^{-1}\mathcal{A}$ with $\mu=1$ and Fig. \ref{Fig2} for $\mu=0.1$.
From the spectral distribution in Figs. \ref{Fig1}-\ref{Fig2}, four preconditioners $P_{SS}$, $P_{RSS}$, $P_{PPSS}$
and $P_{Aug}$ improve the spectral distribution of the original
coefficient matrix $\mathcal{A}$.
As we observe, the eigenvalues of $P_{SS}^{-1}\mathcal{A}$ are $P_{RSS}^{-1}\mathcal{A}$ better clustered than the two other preconditioned matrices.  Moreover, the spectral
distribution of $P_{SS}^{-1}\mathcal{A}$ and $P_{RSS}^{-1}\mathcal{A}$ are almost in line with the theoretical
results, see Theorem \ref{Thm3} and Theorem \ref{Thm4}.

\begin{figure}
    \centering
    \includegraphics[width=6in,height=3.5in]{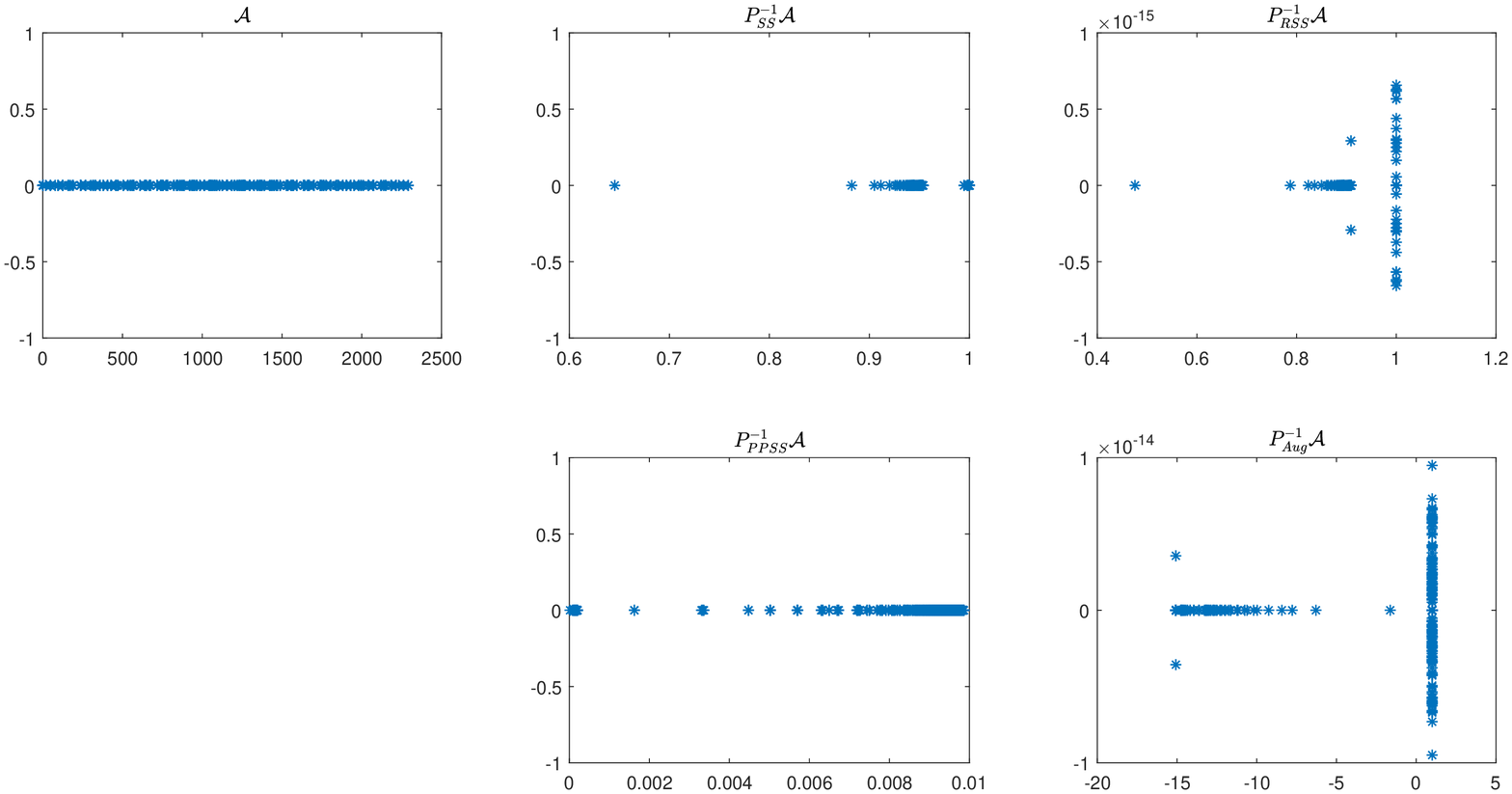}
    \caption{Spectra distribution of Example \ref{Ex1} for $s=16$  with $\mu=1$ and $k=2$.\label{Fig1}}
\end{figure}

\begin{figure}
    \centering
    \includegraphics[width=6in,height=3.5in]{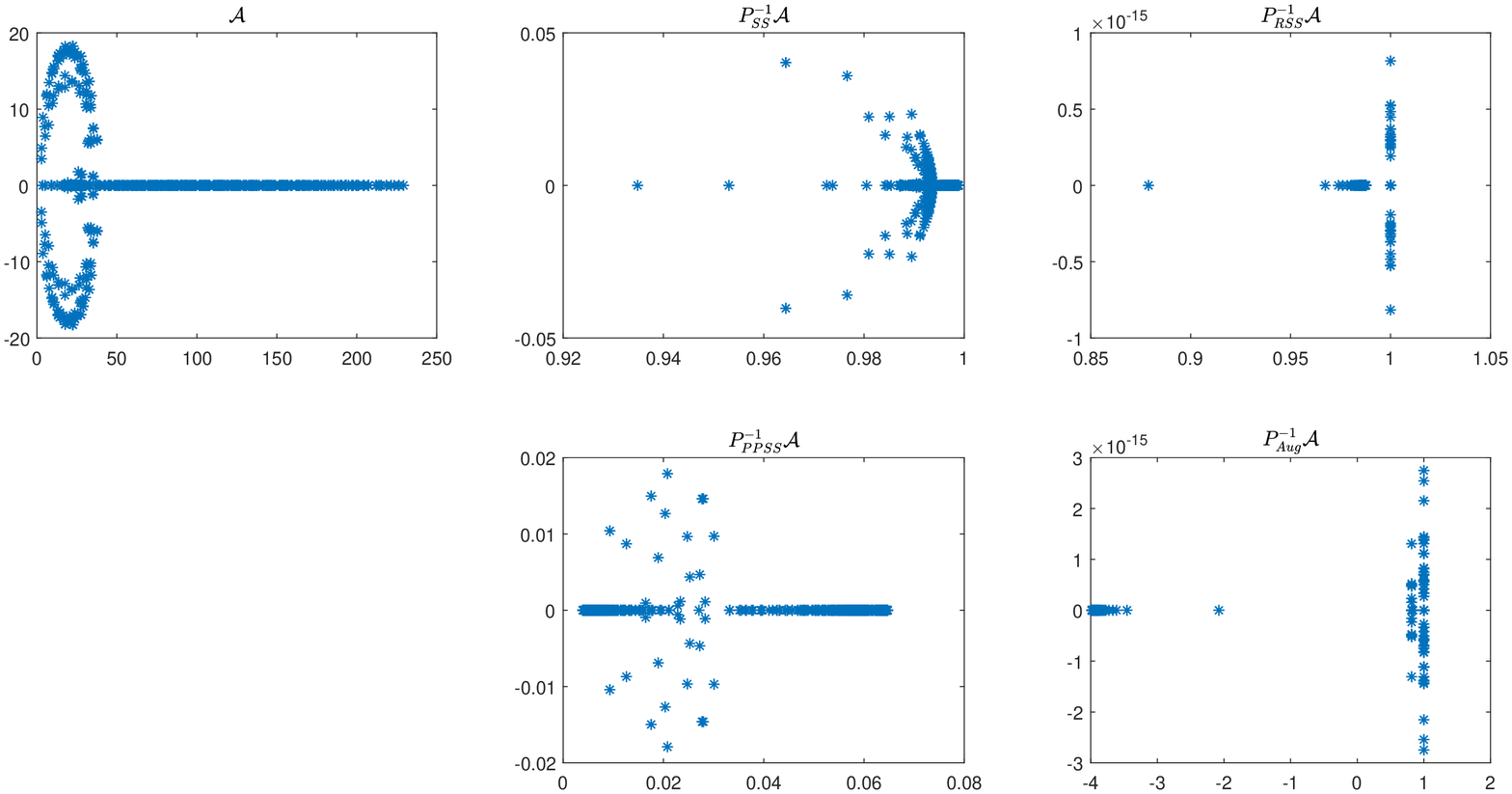}
    \caption{Spectra distribution of Example \ref{Ex1} for $s=16$  with $\mu=0.1$ and $k=2$.\label{Fig2}}
\end{figure}

\end{Example}

\begin{Example}\label{Ex2}\rm
We use the matrix {\bf Szczerba$\slash$Ill\_Stokes} from the UF Sparse Matrix Collection\footnote{https://www.cise.ufl.edu$\slash$research$\slash$sparse$\slash$matrices$\slash$Szczerba$\slash$Ill\_Stokes.html}, which is an ill-conditioned matrix arisen from computational fluid dynamics problems. Generic properties of the test matrix are given in Table \ref{Tbl4}.
The FGMRES (GMRES) method without preconditioning fails to converge in 1000 iterations. So, we present the numerical results of the FGMRES method with the preconditioners $P_{SS}$, $P_{RSS}$, $P_{PPSS}$ and $P_{Aug}$ for different values of the parameter $\alpha$ in Table \ref{Tbl5}. As we observe all the preconditioners reduce the number of iterations of the GMRES method. The minimum value of the CPU time for each of the preconditioner have been underlined. As we see the minimum value of the CPU time is due to the $P_{SS}$ preconditioner. \cred{Numerical results of the preconditioners  $P_{SS}$ and $P_{RSS}$ have been presented in Table \ref{Tbl6}. As we there is a good agreement between the results of the $P_{SS}$ and $P_{RSS}$ preconditioners with $\alpha_*$ and those of with $\alpha_{est}$.} 

\end{Example}
\begin{table}
    \centering\caption{Matrix properties  for Example \ref{Ex2}.\label{Tbl4}}\vspace{-0.5cm}
    \begin{center}
        \label{exact}
        \scalebox{0.8}
        {
            \begin{tabular}{|c|c|c|c|c|c|} \hline
                Matrix                       &  $n$  &  $m$    & $nnz(A)$ & $nnz(B)$ & $nnz(C)$ \\ \hline
                Szczerba$\slash$Ill\_Stokes  & 15672 &  5224   & 73650    & 58242    & 59476    \\ \hline
            \end{tabular}
        }
    \end{center}
\end{table}

\begin{table}
    \centering\caption{Numerical results for Example \ref{Ex2} for different values of $\alpha$.\label{Tbl5}}\vspace{-0.5cm}
    \begin{center}
        \label{exact}
        \scalebox{0.80}
        {
        \begin{tabular}{|c|ccc|ccc|ccc|ccc|} \hline
      & \multicolumn{3}{c|}{$P_{SS}$} & \multicolumn{3}{c|}{$P_{RSS}$} & \multicolumn{3}{c|}{$P_{PPSS}$}  &   \multicolumn{3}{c|}{$P_{Aug}$}\\ \hline
    $\alpha$  & Iters & CPU  & $R_k$ & Iters & CPU  & $R_k$ &  Iters & CPU     & $R_k$ & Iters & CPU     & $R_k$ \\ \hline
    0.1     & 471 & 31.32& 9.8e-8 & 179 & 15.11 & 9.7e-8 & 467   & 22.97   & 1.0e-7  & 188 & 8.31 & 9.9e-8    \\
    0.05    & 358 & 22.14& 9.9e-8 & 169 & 13.57 & 9.8e-8 & 348   & 13.65   & 9.7e-8  & 180 & 7.94 & 9.6e-8   \\
    0.01    & 210 & 13.56& 9.7e-8 & 145 & 11.44 & 9.7e-8 & 164   & 5.19    & 9.1e-8  & 171 & 7.23 & 9.9e-8    \\
    0.005   & 173 & 11.40& 9.6e-8 & 134 & 10.12 & 9.9e-8 & 121   & \underline{{\bf 4.44}}  & 9.2e-8 & 175 & \underline{{\bf 7.19}} & 9.9e-8  \\
    0.001   & 115 & 7.49 & 9.9e-8 & 110 & 7.21  & 9.5e-8 &  66   & 6.97    & 9.1e-8 & 193 & 7.61 & 9.8e-8\\
    0.0005  & 99  & 5.84 & 9.8e-8 & 97  & 5.93  & 9.8e-8 &  66   & 15.23   & 9.0e-8 & 208 & 8.20 & 9.7e-8\\
    0.0001  & 64  & \underline{{\bf 3.91}} & 9.7e-8 & 63 & \underline{{\bf 3.95}}  & 1.0e-7 &  84   & 123.13  & 9.4e-8 & 311 & 16.30 & 1.0e-7\\
   0.00005  & 62  & 4.26 & 9.6e-8 & 61  & 4.23 & 1.0e-7  &  92   & 168.02  &  9.3e-8 & 313 & 18.61 & 9.9e-8\\
   
    \hline
        \end{tabular}
        }
    \end{center}\vspace{0.5cm}
   \cred{ 
    \centering\caption{Numerical results for Example \ref{Ex2} for  $\alpha_{est}$.\label{Tbl6}}\vspace{-0.5cm}
    \begin{center}
    	\label{exact}
    	\scalebox{0.90}
    	{
    		\begin{tabular}{|cccc|cccc|} \hline
    			\multicolumn{4}{|c|}{$P_{SS}$} & \multicolumn{4}{c|}{$P_{RSS}$} \\ \hline
    			$\alpha_{est}$  & Iters & CPU  & $R_k$ & $\alpha_{est}$ & Iters & CPU  & $R_k$ \\ \hline
    			0.000169     & 74 & 4.57 & 9.5e-8 & 0.000169 & 73  & 4.24 & 9.9e-8   \\
    			\hline
    		\end{tabular}
    	}
    \end{center}}

\end{table}

%\begin{table}
%    \centering\caption{Numerical results for Example \ref{Ex2} for  $\alpha_{est}$.\label{Tbl6}}\vspace{-0.5cm}
%    \begin{center}
%        \label{exact}
%        \scalebox{0.90}
%        {
%        \begin{tabular}{|cccc|cccc|} \hline
%       \multicolumn{4}{|c|}{$P_{SS}$} & \multicolumn{4}{c|}{$P_{RSS}$} \\ \hline
%    $\alpha_{est}$  & Iters & CPU  & $R_k$ & $\alpha_{est}$ & Iters & CPU  & $R_k$ \\ \hline
%    0.000169     & 74 & 4.57 & 9.5e-8 & 0.000169 & 73  & 4.24 & 9.9e-8   \\
%    \hline
%        \end{tabular}
%        }
%    \end{center}
%\end{table}

\section{Conclusion}\label{Sec4}

For the asymmetric saddle point problems, we have
presented the shift-splitting preconditioner and its relaxed version
to improve the convergence speed of Krylov subspace method (such as
GMRES/FGMRES). The eigenvalue distribution of the related preconditioned
matrices have been provided. Moreover, we have proved that the shift-splitting iteration method for
the asymmetric saddle point problems
is convergent under suitable conditions. Numerical experiments from
the Stokes problem are given to verify the efficiency of the
shift-splitting preconditioner and its relaxed version.

\section*{ Acknowledgment}

\cred{The authors would like to thank the anonymous referee for helpful comments and suggestions.}

 {
}
\end{document}